# THE ARCTIC CIRCLE BOUNDARY AND THE AIRY PROCESS

By Kurt Johansson[1]

*Royal Institute of Technology*

We prove that the, appropriately rescaled, boundary of the north polar region in the Aztec diamond converges to the Airy process. The proof uses certain determinantal point processes given by the extended Krawtchouk kernel. We also prove a version of Propp's conjecture concerning the structure of the tiling at the center of the Aztec diamond.

**1. Introduction and results.** Domino tilings of the Aztec diamond were introduced in [8, 9]. Asymptotic properties of random domino tilings of the Aztec diamond have been studied in [5, 12, 15]. In particular, in [12] the existence of the so-called arctic circle was proved. The arctic circle is the asymptotic boundary of the disordered so-called temperate region of the tiling. Outside this boundary the tiling forms a completely regular brick wall pattern. The methods in [12] combined with the results in [13] show that the fluctuations of the point of intersection of the boundary of the temperate region with a line converge to the Tracy–Widom distribution of random matrix theory. In this paper we extend this result to show that the fluctuations of the boundary around the arctic circle converges to the Airy process introduced in [23]. The paper is a continuation of the approach used in [14] and [15], where certain point processes with determinantal correlation functions [24] and the Krawtchouk ensemble, were used. We will use the general techniques developed in [16] and investigate an extended point process which also has determinantal correlation functions given by a kernel, which we call the extended Krawtchouk kernel.

The *Aztec diamond*, $A_n$, of order $n$ is the union of all lattice squares $[m, m+1] \times [l, l+1]$, $m, l \in \mathbb{Z}$, that lie inside the region $\{(x_1, y_1); |x_1| + |y_1| \leq$

Received June 2003; revised April 2004.

[1]Supported in part by the Swedish Science Research Council and the Göran Gustafsson Foundation (KVA).

*AMS 2000 subject classifications.* Primary 60K35; secondary 82B20, 15A52.

*Key words and phrases.* Airy process, determinantal process, Dimer model, random matrices, random tiling.

$n+1\}$. A *domino* is a closed $1 \times 2$ or $2 \times 1$ rectangle in $\mathbb{R}^2$ with corners in $\mathbb{Z}^2$, and a *tiling* of a region $R \subseteq \mathbb{R}^2$ by dominoes is a set of dominoes whose interiors are disjoint and whose union is $R$. Let $\mathcal{T}(A_n)$ denote the set of all domino tilings of the Aztec diamond. The basic coordinate system used here will be referred to as coordinate system I (CS-I).

Color the Aztec diamond in a checkerboard fashion so that the leftmost square in each row in the top half is white. A horizontal domino is *north-going* (N) if its leftmost square is white, otherwise it is *south-going* (S). Similarly, a vertical domino is *west-going* (W) if its upper square is white, otherwise it is *east-going* (E). Two dominoes are *adjacent* if they share an edge, and a domino is adjacent to the boundary if it shares an edge with the boundary of the Aztec diamond. The *north polar region* (NPR) is defined to be the union of those north-going dominoes that are connected to the boundary by a sequence of adjacent north-going dominoes. The south, west and east polar regions are defined analogously. In this way a domino tiling partitions the Aztec diamond into four polar regions, where we have a regular brick wall pattern, and a fifth central region, the *temperate region*, where the tiling pattern is irregular.

Let $T \in \mathcal{T}(A_n)$ be a tiling of the Aztec diamond and let $v(T)$ denote the number of vertical dominoes in $T$. We can define a probability measure on $\mathcal{T}(A_n)$ by letting vertical dominos have weight $a$ and horizontal dominos weight 1, that is,

$$(1.1) \qquad \mathbb{P}[T] = \frac{a^{v(T)}}{\sum_{T \in \mathcal{T}(A_n)} a^{v(T)}}.$$

In this paper we will prove the asymptotic results for the uniform case ($a = 1$) only to keep the asymptotic analysis simpler.

We will study the part of the boundary of the NPR which lies above a neighborhood of $x_1 = 0$. To define the boundary we will use certain nonintersecting paths which describe the domino configuration and which are also essential in our analysis below. They were called DR-paths ([25], page 277) in [15]. On a W-domino placed so that it has corners at $(0,0)$ and $(1,2)$, we draw a line from $(0,1/2)$ to $(1,3/2)$. On an E-domino placed in the same position, we draw a line from $(0,3/2)$ to $(1,1/2)$, and on an S-domino placed so that it has corners at $(0,0)$ and $(2,1)$, we draw a line from $(0,1/2)$ to $(2,1/2)$. We do not draw any line on an N-domino. As discussed in [15], these lines will form nonintersecting paths from $A_r = (-n-1+r, -r+1/2)$ to $B_r = (n+1-r, -r+1/2)$, $1 \le r \le n$. The top curve, from $A_1$ to $B_1$, can be viewed as a function, $t \to X_n(t)$, $|t| \le n$, in CS-I. We will call $X_n(t)$ the NPR-boundary process, see Figure 1. The NPR is exactly the part of the domino tiling that lies above $X_n(t)$, and consists only of N-dominoes. Between the nonintersecting paths there are other regions of N-dominoes.

Before we formulate the limit theorem for $X_n(t)$, we recall the definition of the Airy process. The *extended Airy kernel* is defined by

$$
(1.2) \quad \begin{aligned} &A(\tau, \xi; \tau', \xi') \\ &= \begin{cases} \int_0^\infty e^{-\lambda(\tau-\tau')} \operatorname{Ai}(\xi+\lambda) \operatorname{Ai}(\xi'+\lambda) \, d\lambda, & \text{if } \tau \geq \tau', \\ -\int_{-\infty}^0 e^{-\lambda(\tau-\tau')} \operatorname{Ai}(\xi+\lambda) \operatorname{Ai}(\xi'+\lambda) \, d\lambda, & \text{if } \tau < \tau', \end{cases} \end{aligned}
$$

for $\tau, \tau', \xi, \xi' \in \mathbb{R}$. Fix $\tau_1 < \cdots < \tau_m$ and let $\Lambda_m = \{\tau_1, \ldots, \tau_m\}$. Let $\mu$ be the product of counting measure on $\Lambda_n$ and Lebesgue measure on $\mathbb{R}$. Define $f$ on $\Lambda_m \times \mathbb{R}$ by $f(\tau_j, x) = \mathbb{1}_{(\xi_j, \infty)}(x)$ for given numbers $\xi_1, \ldots, \xi_m$. It is proved in [16] that $f^{1/2}(\tau, x) A(\tau, x; \sigma, y) f^{1/2}(\sigma, y)$ is the integral kernel of a trace class operator on $L^2(\Lambda_m \times \mathbb{R}, \mu)$. The *Airy process*, $t \to \mathcal{A}(t)$, is the stationary stochastic process whose finite-dimensional distributions are given by

$$(1.3) \quad \mathbb{P}[\mathcal{A}(\tau_1) \leq \xi_1, \ldots, \mathcal{A}(\tau_m) \leq \xi_m] = \det(I - f^{1/2} A f^{1/2})_{L^2(\Lambda_m \times \mathbb{R})}.$$

It is proved in [23], see also [16], that $\mathcal{A}(\tau)$ has continuous paths. The distribution of $\mathcal{A}(\tau)$ is $F_2$, the Tracy–Widom distribution for the largest eigenvalue of a GUE-matrix [27]. It has recently been shown that the distribution function in (1.3) satisfies certain differential equations [1, 28].

Our main result is the following:

THEOREM 1.1. *Let $X_n(t)$ be the NPR-boundary process and $\mathcal{A}(\tau)$ the Airy process as defined above and let the weight $a$ in (1.1) be equal to $1$. Then,*

$$(1.4) \quad \frac{X_n(2^{-1/6} n^{2/3} t) - n/\sqrt{2}}{2^{-5/6} n^{1/3}} \to \mathcal{A}(t) - t^2,$$

*as $n \to \infty$, in the sense of convergence of finite-dimensional distributions.*

The theorem could be extended to a general weight $a$ in (1.1) by the same type of argument. We restrict to $a = 1$ for simplicity. Similarly, it is possible to show convergence along other parts of the boundary, except right near the point where the arctic circle is tangent to the asymptotic square containing the tiling where the boundary behavior is different. See the remark in the last paragraph of Section 2.

As was done in [16] for the convergence of the interface of a polynuclear growth model to the Airy process, it should be possible to extend Theorem 1.1 to a functional limit theorem. It would then follow that $\max_t X_n(t)$, suitably rescaled, asymptotically has $F_1$ fluctuations, where $F_1$ is the Tracy–Widom distribution for the largest eigenvalue of a GOE-matrix [27].

The above result also has an interpretation as a convergence theorem for a certain polynuclear growth model, see [15], Section 2.4. The polynuclear growth model studied in [23] is a limiting version of this model ($a \to 0$, $n \to \infty$ at appropriate tates). The NPR-boundary is also directly related to the shape in the corner growth model studied in [13] and [15], Section 2.4, and, hence, also to the totally asymmetric exclusion process. If $G(M, N)$ is the last passage time as in [13], the results of the present paper show that the fluctuations of the boundary of the shape $\Omega_t = \{(M, N); G(M, N) + M + N - 1 \leq t\}$ close to the diagonal $(N, N)$, when $q = 1/2$, converges to the Airy process. This can be extended to the other parts of the boundary, away from the axes, and to other values of $q$. This extends the relation between the Meixner ensemble and the Krawtchouk ensemble given in Lemma 2.9 in [15], see also [20].

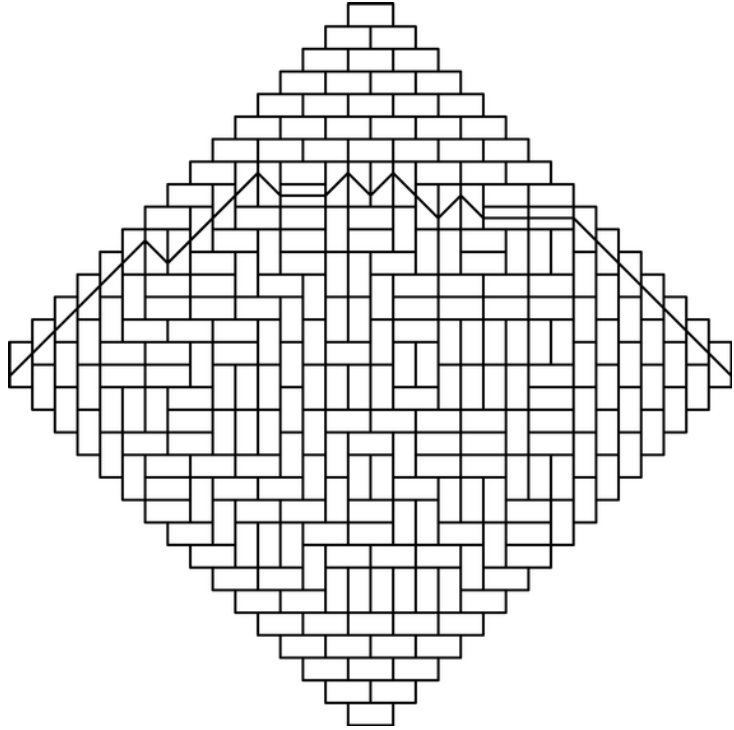

Fig. 1. *An NPR-boundary process.*

The theorem will be proved using a certain determinantal point process. One way of seeing this point process is to put a *green particle* at the center of the black square on each S-domino and W-domino, and a *red particle* at the center of the white square of each S- and W-domino. These dots define a point process and we will see below that it has determinantal correlation functions. It is directly related to the nonintersecting paths defined above and also to the zig-zag paths around black and white squares discussed in [8, 9, 14]. The precise definition of this point process will be given in Section

2. The asymptotic results needed to prove Theorem 1.1 will be discussed in Section 3. It has been conjectured in [12] that at the center of a uniform random tiling of the Aztec diamond the tiling looks like a random domino tiling of the plane under the Burton–Pemantle measure [4]. A proof of a version of this conjecture will be given in Section 4.

It is natural to conjecture that boundaries between irregular and regular tiling regions in many other two-dimensional random tiling problems are also described by the Airy process. In terms of the so-called height function, which we will not define here, this would mean that we would see the Airy process where we have a boundary between a flat surface and a curved surface. This type of result can also be proved in the rhombus tiling problem discussed in [22]. This has been done recently in [10]. In this problem it would also be possible to use the formulas derived in [22], which would lead to computations very similar to those in the present paper. Another candidate where it may be possible to prove convergence to the Airy process would be for fluctuations around the arctic circle (ellipse) in rhombus tilings of hexagons [7, 15], but here the asymptotic analysis appears to be more difficult [2].

**2. The point process.** Introduce a new coordinate system (CS-II) with origin at $(-n, -1/2)$ and axes $e_{II} = (1,1)$, $f_{II} = (-1,1)$, which gives the following coordinate transformation between CS-I and CS-II,

$$\begin{aligned} x_1 &= x_2 - y_2 - n, \\ y_1 &= x_2 + y_2 - 1/2. \end{aligned} \quad (2.1)$$

In CS-II the nonintersecting paths defined in Section 1 go from $A_j = (0, -j+1)$ to $B_j = (n+1-j, -n)$, $1 \leq j \leq n$, and have three types of steps $(1,0)$, $(0,-1)$ and $(1,-1)$, see Figure 2. These nonintersecting paths specify the tiling uniquely. The measure (1.1) is obtained by letting the steps $(1,0)$ and $(0,-1)$, which correspond to vertical dominos, have weight $a$ and the step $(1,-1)$ weight 1.

To formalize this, let $\mathcal{G} = (V, E)$ be a directed graph with vertex set $V = \mathbb{N} \times \mathbb{Z}$ and directed edges from $(i,j)$ to $(i+1,j)$, $(i, j-1)$ and $(i+1, j-1)$ for $i \geq 1$, and from $(0,j)$ to $(1, j-1)$ and $(1,j)$. The edges from $(i,j)$ to $(i+1, j-1)$ have weight 1, whereas all other edges have weight $a$. A path $\pi$ from $A \in V$ to $B \in V$ is a directed path along succesive directed edges starting at $A$ and ending at $B$. The weight $\omega(\pi)$ of a path $\pi$ is the product of the weights of all the edges in the path. Two paths $\pi$ and $\pi'$ are nonintersecting if they do not share a common vertex. If $\underline{A} = (A_1, \ldots, A_m)$ and $\underline{B} = (B_1, \ldots, B_m)$, where $A_r, B_r \in V$, then $\mathcal{P}_{\text{n.i.}}(\underline{A}; \underline{B})$ denotes the set of all nonintersecting paths $\pi_1, \ldots, \pi_m$, where $\pi_j$ goes from $A_j$ to $B_j$. The weight of a family $\underline{\pi} = (\pi_1, \ldots, \pi_m)$ of paths is $\omega(\underline{\pi}) = \prod_{j=1}^m \omega(\pi_j)$. Take $N \geq n$ and set $A_j = (0, 1-j)$ and $C_j = (n, -n+1-j)$, $1 \leq j \leq N$; see Figure 3.

LEMMA 2.1. *If $(\pi, \ldots, \pi_N)$ are paths in $\mathcal{P}_{\mathrm{n.i.}}(\underline{A}; \underline{C})$, then $\pi_k$ goes through the vertex $B_k = (n+1-k, -n)$, $1 \le k \le n$.*

PROOF. Since $\pi_1$ ends at $C_1 = B_1$, the claim is true for $k = 1$. Each $\pi_k$, $1 < k \le n$, has to go through one of the points $B-2, \ldots, B_n$ since $C_2, \ldots, C_n$ lie below the line $y_2 = -n$ and all the paths are right/down directed. The nonintersection constraint implies that $\pi_{k+1}$ has to pass this line to the left of $\pi_k$ and the claim follows. □

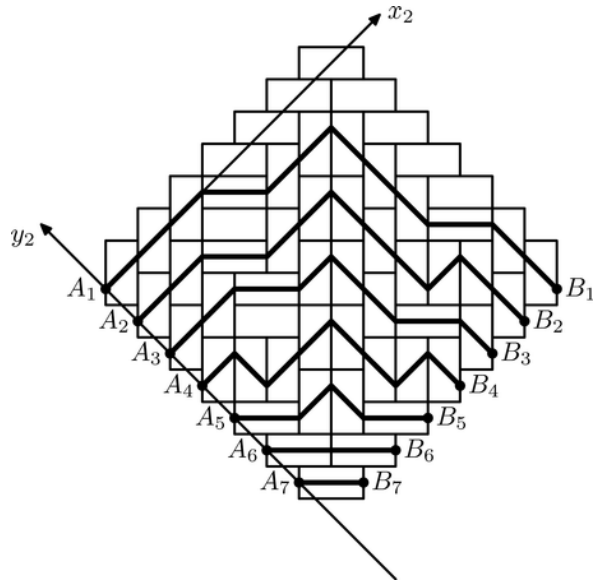

FIG. 2. *CS-II and nonintersecting paths describing the tiling.*

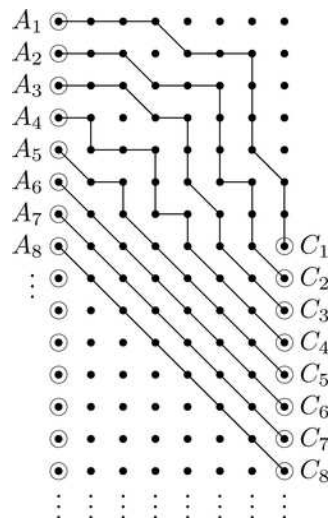

FIG. 3. *The nonintersecting paths in the graph $\mathcal{G}$.*

By this lemma there are well-defined projections

$$P : \Omega \doteq \mathcal{P}_{\text{n.i.}}(A_1, \ldots, A_N; C_1, \ldots, C_N) \to \Omega' \doteq \mathcal{P}_{\text{n.i.}}(A_1, \ldots, A_n; B_1, \ldots, B_n)$$

and

$$Q : \Omega \to \Omega'' \doteq \mathcal{P}_{\text{n.i.}}(B_1, \ldots, B_n, A_{n+1}, \ldots, A_N; C_1, \ldots, C_N),$$

such that $\omega(\pi) = \omega(P(\pi))\omega(Q(\pi))$. An event $D$ in the domino tiling of the size $n$ Aztec diamond corresponds, via the bijection, to an event $\widetilde{D} \subseteq \Omega'$. Set $D_N^* = \{\pi \in \Omega; P(\pi) \in \widetilde{D}\}$. Then,

$$\begin{aligned}
(2.2) \quad \mathcal{P}_N[D_N^*] &= \frac{\sum_{\underline{\pi} \in D_N^*} \omega(\underline{\pi})}{\sum_{\underline{\pi} \in \Omega} \omega(\underline{\pi})} \\
&= \frac{(\sum_{\underline{\pi}' \in \widetilde{D}} \omega(\underline{\pi}'))(\sum_{\underline{\pi}'' \in \Omega''} \omega(\underline{\pi}'))}{(\sum_{\underline{\pi}' \in \Omega'} \omega(\underline{\pi}''))(\sum_{\underline{\pi}'' \in \Omega''} \omega(\underline{\pi}''))} = \mathbb{P}[D].
\end{aligned}$$

The right-hand side in (2.2) is independent of $N \geq n$ and, hence,

$$(2.3) \qquad \mathbb{P}[D] = \lim_{N \to \infty} \mathbb{P}_N^*[D_N^*].$$

We want to map, bijectively and preserving weights, the nonintersecting paths in $\Omega$ to a new family of nonintersecting paths which is appropriate for the application of the general results in [16]. The idea is as follows. Each path $\pi_k$ from $A_k$ to $C_k$ has a first and a last vertex on each vertical line $x_2 = j$. The *first* and the *last vertex* could, of course, coincide. We want to put the first and the last vertex on different vertical lines. These first and last vertices form the point process we are interested in. Let $\mathcal{G}' = (V', E')$ be the directed graph with $V' = \mathbb{N} \times \mathbb{Z}$ and directed edges from $(2i-1, j)$ to $(2i, j)$, $i \geq 1$, from $(2i, j)$ to $(2i, j-1)$, from $(2i, j)$ to $(2i+1, j)$, from $(2i, j)$ to $(2i+1, j+1)$, $i \geq 1$, from $(0, j)$ to $(1, j)$ and from $(0, j)$ to $(1, j+1)$, $j \in \mathbb{Z}$. We put the weight $a$ on the edges from $(2i, j)$ to $(2i+1, j+1)$ and from $(2i, j)$ to $(2i, j-1)$, $i \geq 1$, $j \in \mathbb{Z}$. All other edges have weight $= 1$.

We can describe a path $\pi_k$ from $A_k$ to $C_k$ by giving the first, $P_k(j)$, and the last, $Q_k(j)$, point on each verical line $x_2 = j$, $0 \leq j \leq n$; $P_k(0) = Q_k(0)$ always. From $Q_k(j)$ to $P_k(j+1)$ we take either a step $(1, 0)$ or a step $(1, -1)$, and from $P_k(j)$ to $Q_k(j)$ we take a certain number, $\geq 0$, of down steps $(0, -1)$. Map $Q_k(j) = (j, q)$ to $R_k(2j) = (2j, q+j)$, $0 \leq j \leq n$, and $P_k(j) = (j, p)$ to $R_k(2j-1) = (2j-1, p+j)$. A step from $Q_k(j)$ to $P_k(j+1)$ is mapped to two steps. One step from $R_k(2j)$ to $R_k(2j+1)$, which is either $(1, 0)$ or $(1, 1)$, and then a fixed step from $R_k(2j+1) = (2j+1, p+j)$ to $(2j+2, p+j)$, that is, a step $(1,0)$. The vertical steps from $P_k(j)$ to $Q_k(j)$, if there are any, are mapped to the same number of vertical steps from $(2j, p+j)$ to $(2j, q+j)$. Set $A'_k = (2n, 1-k)$, $1 \le k \le N$. The above procedure maps $\mathcal{P}_{\text{n.i.}}(A_1, \ldots, A_N; C_1, \ldots, C_N)$ to $\mathcal{P}_{\text{n.i.}}(A_1, \ldots, A_N; A'_1, \ldots, A'_N)$ bijectively, see Figure 4. Also, the weight of a family of paths is preserved. The probability of a family of paths $\underline{\pi}$ in $\widehat{\Omega}_N \doteq \mathcal{P}_{\text{n.i.}}(\underline{A}; \underline{A}')$ is

$$\widehat{\mathbb{P}}_N(\underline{\pi}) = \frac{\omega(\underline{\pi})}{\sum_{\underline{\pi} \in \widehat{\Omega}_N} \omega(\underline{\pi})}. \tag{2.4}$$

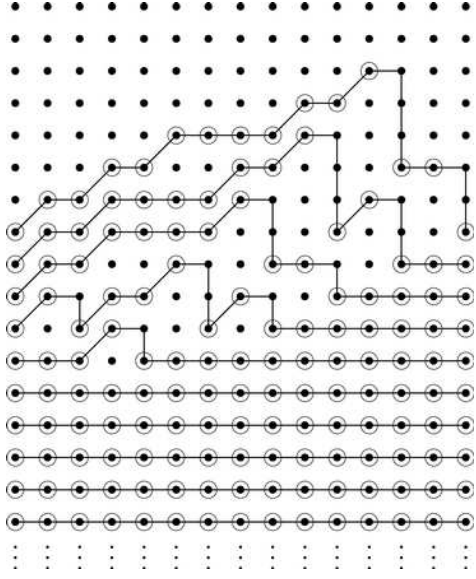

FIG. 4. *The nonintersecting paths in the graph $\mathcal{G}'$ corresponding to the tiling in Figure 2. The particles in the determinantal process are the circled dots.*

An event $D_N^*$ in $\Omega$ as above maps to an event $\widehat{D}_N$ in $\widehat{\Omega}_N$ using the bijection, and from (2.3), we have $\mathbb{P}[D] = \lim_{N \to \infty} \widehat{\mathbb{P}}_N(\widehat{D}_N)$.

The paths $\pi_1, \ldots, \pi_N$ are uniquely specified by the points $(R_k(j))_{1 \le j \le 2n-1, 1 \le k \le N}$. We will write $R_k(j) = (j, x_k^j)$, $1 \le j \le 2n-1$, $1 \le k \le N$ and $\underline{x} = (x_k^j)$; $x_k^0 = x_k^{2n} = 1-k$, $k \ge 1$; $x_1^j > x_2^j > \cdots > x_N^j$. These points define a particle configuration in $\{1, \ldots, 2n-1\} \times \mathbb{Z}$ and we obtain a point process on this space which we can think of as $2n-1$ copies of $\mathbb{Z}$. We see that the NPR-boundary process is obtained by joining $Q_1(0), P_1(0), P_1(1), Q_1(1), \ldots, P_1(n), Q_1(n)$ with straight lines in this order. If $Q_1(j)$ has coordinates $(j, q_j)$ in CS-II, then $R_1(2j) = (2j, q_j + j)$, so $q_j = x_1^{2j} - j$. Similarly, if $P_1(j)$ has coordinates

$(j, p_j)$ in CS-II, then $p_j = x_1^{2j-1} - j$. We see from (2.1) that

(2.5)
$$Q_1(j) = (2j - x_1^{2j} - n, x_1^{2j} - 1/2),$$
$$P_1(j) = (2j - x_1^{2j-1} - n, x_1^{2j-1} - 1/2)$$

in CS-I.

The red points described in the introduction are $P_k(j) + (-1/2, 0)$ in CS-I and the green points are $Q_k(j) + (1/2, 0)$; note that we can have $Q_k(j) = P_k(j)$. This is seen from how we defined the nonintersecting paths and the red and green points in Section 1.

We can define *zig-zag paths* around white and black squares as in [8, 9, 14, 15]. Consider the sequence of white squares in the Aztec diamond with opposite corners $Q_k^r = (-r + k, n + 1 - k - r)$, $k = 0, \ldots, n+1$, where $r$, $1 \leq r \leq n$, is fixed. The zig-zag steps from $Q_k^r$ to $Q_{k+1}^r$ go either one unit step to the right and then one unit step down (ES-step) or the other way around (SE-step), and in such a way that it does not intersect a domino. Similarly, we can define zig-zag paths around black squares between $(-r, n - r)$ and $(n - r, -r)$. It is not hard to see, compare the proof of Lemma 2.2 in [15], that the zig-zag paths around white squares have a red dot exactly when we have an ES-step around it, and that zig-zag paths around black squares have SE-steps around squares with green dots. In this way each zig-zag path corresponds to a unique particle configuration, and we can view the point process as the totality of all zig-zag particles.

It was proved in [14] and [15] that the zig-zag particles along a single line define a point process described by the Krawtchouk ensemble. The possible positions of the red particles can be taken to be $\{0, 1, \ldots, n\}$ and the $r$th zig-zag path has $r$ red particles at $h_1, \ldots, h_r$ (no ordering). The probability of having particles at exactly these points is [15], Theorem 2.2,

(2.6)
$$\mathbb{P}[h] = \frac{1}{Z_{r,n,q}} \Delta_r(h)^2 \prod_{j=1}^r \binom{n}{h_j} q^{h_j} p^{n-h_j},$$

where $q = a^2(1 + a^2)^{-1}$, $p = 1 - q$ and $\Delta_r(h)$ is the Vandermonde determinant. This *Krawtchouk ensemble* has determinantal correlation functions given by the *Krawtchouk kernel*,

(2.7) $$K_{r,n,q}(x, y) = \sum_{k=0}^{r-1} p_k(x; q, n) p_k(y; , q, n) \left[ \binom{n}{x} q^x p^{n-x} \binom{n}{y} q^y p^{n-y} \right]^{1/2},$$

where

(2.8) $$p_k(x; q, n) = \binom{n}{k}^{-1/2} (qp)^{-k/2} (-1)^k \frac{1}{2\pi i} \int_\gamma \frac{(1 + pz)^x (1 - qz)^{n-x}}{z^k} \frac{dz}{z}$$

are multiples of the ordinary Krawtchouk polynomials; $\gamma$ is a circle centered at the origin with radius $\leq \min(1/p, 1/q)$. $p_k(x; q, n) = \binom{n}{k}^{1/2} (q/p)^{k/2} K_k(x; q, n)$, where $K_k$ are the standard Krawtchouk polynomials, see, for example, [19].

Next, we want to show that the fact that we have a determinantal point process when we restrict to a single line can be extended to show that the whole point process is determinantal. Using the Lindström–Gessel–Viennot method, the combinatorial version of the Karlin–McGregor theorem (see, e.g., [26]) we can write the probability for a certain particle configuration $\underline{x}$ as a product of determinants

$$(2.9) \quad p_{N,n}(\underline{x}) = \frac{1}{Z_{N,n}} \prod_{r=0}^{2n-1} \det(\phi_{r,r+1}(x_j^r, x_k^{r+1}))_{j,k=1}^N,$$

where the transition function $\phi_{r,r+1}(x, y)$ gives the weight of all paths going from $x$ on the vertical line $r$ to $y$ on the vertical line $r + 1$. The right-hand side of (2.9) is a symmetric function of $x_1^r, \ldots, x_N^r$ for any $r$, so we need not consider the variables on a single vertical line as ordered. From the definitions of the directed graph $\mathcal{G}'$ and its weights, we see that we can describe the nonintersecting paths using two types of transitions between different vertical lines. From vertical line $2i$ to $2i + 1$, we can go from $(2i, j)$ to $(2i + 1, j)$, with weight 1, or from $(2i, j)$ to $(2i + 1, j + 1)$, with weight $a$. From line $2i + 1$ to $2i + 2$, we have to take a step to the right, with weight 1, and then a certain number $\geq 0$ of down steps, with weight $a$. This leads to the following transition functions:

$$(2.10) \quad \phi_{2i,2i+1}(x, y) = \alpha(y - x) = \begin{cases} a, & \text{if } y - x = 1, \\ 1, & \text{if } y - x = 0, \\ 0, & \text{otherwise,} \end{cases}$$

and

$$(2.11) \quad \phi_{2i+1,2i+2}(x, y) = \beta(y - x) = \begin{cases} a^{-(y-x)}, & \text{if } y - x \leq 0, \\ 0, & \text{otherwise.} \end{cases}$$

Set $f_{2i}(z) = az + 1$ and $f_{2i+1}(z) = (1 - a/z)^{-1}$ so that $\alpha(n) = \hat{f}_{2i}(n)$ and $\beta(n) = \hat{f}_{2i+1}(n)$, where $\hat{f}_r(n)$ is the $n$th Fourier coefficient of $f_r$. We assume that $0 < a < 1$. The case $a = 1$ will be handled by taking a limit $a \to 1-$. Note that $f_r$ has a Wiener–Hopf factorization $f_r = f_r^+ f_r^-$, where $f_{2i}^+(z) = az + 1$, $f_{2i}^- = f_{2i+1}^+ = 1$ and $f_{2i+1}^-(z) = (1 - a/z)^{-1}$. It follows from [16], Theorem 17, that the point process on $V_n = \{1, \ldots, 2n - 1\} \times \mathbb{Z}$ defined by (2.9) has determinantal correlation functions. The probability of finding particles at $z_j = (r_j, y_j)$, $1 \leq j \leq m$, where $0 < r_j < 2n$ and $y_j \in \mathbb{Z}$, is

$$(2.12) \quad \det(K_{N,n}(z_j; z_k))_{j,k=1}^m,$$

where the kernel is given by the formula

$$(2.13) \quad K_{N,n}(r, x; s, y) = \widetilde{K}_{N,n}(r, x; s, y) - \phi_{r,s}(x, y).$$

Here $\phi_{r,s} = \phi_{r,r+1} * \cdots * \phi_{s-1,s}$ if $s > r$ and $\phi_{r,s} \equiv 0$ if $r \geq s$. Furthermore,

$$(2.14) \quad \widetilde{K}_{N,n}(r,x;s,y) = \sum_{i,j=1}^{N} \phi_{r,2n}(x, x_i^{2n})(A^{-1})_{ij}\phi_{0,s}(x_j^0, y),$$

where $A = (a_{ij})_{i,j=1}^{N}$, $a_{ij} = \phi_{0,2n}(x_i^0, x_j^{2n})$. This is a consequence of the following formula, Proposition 2.1 in [16]. Assume that $g: V_n \to \mathbb{R}$ is a bounded function. Then,

$$(2.15) \quad \mathbb{E}_{N,n}\left[\prod_{r=1}^{2n-1}\prod_{j=1}^{N}(1 + g(r, x_j^r))\right] = \det(I + gK_{N,n})_{L^2(V_n, \mu)},$$

where $\mu$ is counting measure on $V_n$. Here $\mathbb{E}_{N,n}$ is expectation with respect to the probability measure (2.9). If $g(r, x) = 0$, for $x < -n + [r/2]$, then the left-hand side of (2.15) depends only on the part of the particle configuration that corresponds to the Aztec diamond. We know that this is independent of $N$, compare (2.3), and we can replace the kernel $K_{N,n}$ by its limit as $N \to \infty$. This limit is given by, Proposition 1.8 in [16],

$$(2.16) \quad K_n(r,x;s,y) = \widetilde{K}_n(r,x;s,y) - \phi_{r,s}(x,y),$$

where

$$(2.17) \quad \widetilde{K}_n(r,x;s,y) = \frac{1}{(2\pi i)^2} \int_{|z|=r_1} \frac{dz}{z} \int_{|w|=r_2} \frac{dw}{w} \frac{w^y}{z^x} G_{n,r,s}(z,w) \frac{z}{z-w},$$

$a < r_1 < 1/a$, $0 < r_2 < r_1$,

$$(2.18) \quad G_{n,r,s}(z,w) = \frac{\prod_{t=r}^{2n-1} f_t^-(1/z) \prod_{t=0}^{s-1} f_t^+(1/w)}{\prod_{t=0}^{r-1} f_t^+(1/z) \prod_{t=s}^{2n-1} f_t^-(1/w)}$$

and

$$(2.19) \quad \phi_{r,s}(x,y) = \frac{1}{2\pi} \int_{-\pi}^{\pi} e^{i(y-x)\theta} G_{n,r,s}(e^{i\theta}, e^{i\theta})\,d\theta$$

for $r < s$ and $\phi_{r,s}(x,y) = 0$ otherwise.

Inserting the formulas for $f_t^\pm$ into (2.18) we see that, with $\varepsilon_1, \varepsilon_2 \in \{0,1\}$,

$$(2.20) \quad G_{n,2r-\varepsilon_1, 2s-\varepsilon_2}(z,w) = \frac{(1-aw)^{n-s+\varepsilon_2}(1+a/w)^s}{(1-az)^{n-r+\varepsilon_1}(1+a/z)^r}.$$

Using Cauchy's theorem, we can deform the integration contours to $\gamma$ given by $w(t) = \alpha_1 e^{it}$, $-\pi \leq t \leq \pi$, and $\Gamma$ given by $z(t) = \alpha_2 + it$, $t \in \mathbb{R}$, where $0 < \alpha_1 < \alpha_2 < 1/a$. We can then let $a \to 1-$ to get the case $a = 1$ using the continuity of all expressions involved. Then,

$$(2.21) \quad \begin{aligned}\widetilde{K}_n(2r,x;,2s,y) \\ = \frac{1}{(2\pi i)^2} \int_\Gamma \frac{dz}{z} \int_\gamma \frac{dw}{w} \frac{w^y(1-w)^{n-s}(1+1/w)^s}{z^x(1-z)^{n-r}(1+1/z)^r} \frac{z}{z-w}.\end{aligned}$$

It follows that the probability of finding particles at positions $z_j = (r_j, y_j)$, where $y_j \geq -n + [r_j/2]$, $1 \leq j \leq m$, is given by (2.12) with $K_n$ instead of $K_{N,n}$. Consequently, for all our computations we will use $K_n$ given by (2.16). We will call this kernel the *extended Krawtchouk kernel*. For future use we record the following consequence of (2.14). Let $r_1 < \cdots < r_m$ and let $\ell_1, \ldots, \ell_m$ be given real numbers such that $\ell_k \geq -n + [r_k/2]$. Then, for $N \geq n$,

$$
\begin{aligned}
&\mathbb{P}_{N,n}\left[\max_{1 \leq k \leq N} x_k^{r_j} \leq \ell_j, 1 \leq j \leq m\right] \\
&\quad (2.22) \\
&= \mathbb{E}_{N,n}\left[\prod_{r=1}^{2n-1} \prod_{j=1}^{N}(1 + g(r, x_j^r))\right] = \det(I + gK_n)_{L^2(\Lambda_m, \mu)},
\end{aligned}
$$

where $g(r, x) = -\mathbb{1}_{(\ell_k, \infty)}(x)$ if $r = r_k$, $1 \leq k \leq m$, and $g(r, x) = 0$ otherwise. Since $x_1^r = \max_{1 \leq k \leq N} x_k^r$ when the variables are ordered, we can combine (2.22) with (2.5) to study the boundary of the NPR. To use it to prove Theorem 1.1 we need some asymptotic results for the extended Krawtchouk kernel. These will be proved in the next section.

To justify calling (2.16) the extended Krawtchouk kernel we should show how it is related to the Krawtchouk kernel defined by (2.7). Let $p_k$ be the normalized Krawtchouk polynomials given by (2.8) and set $p_k \equiv 0$ if $k < 0$ or $k > n$. Define, $0 < q < 1$, $p = 1 - q$,

$$
(2.23) \quad L_{n,q}(r, x; s, y) = \sum_{k=-\infty}^{-1} \binom{n}{k+s}^{-1/2} \binom{n}{k+r}^{1/2} \\
\times p_{k+r}(x; q, n) p_{k+s}(y; q, n) [w_q(x) w_q(y)]^{1/2}
$$

if $r \leq s$, and

$$
(2.24) \quad L_{n,q}(r, x; s, y) = -\sum_{k=0}^{\infty} \binom{n}{k+s}^{-1/2} \binom{n}{k+r}^{1/2} \\
\times p_{k+r}(x; q, n) p_{k+s}(y; q, n) [w_q(x) w_q(y)]^{1/2}
$$

if $r > s$, where $w_q(x) = \binom{n}{x} q^x p^{y-x}$. Note that $L_{n,q}(r, x; r, y) = K_{r,n,q}(x, y)$, where $K_{r,n,q}$ is the Krawtchouk kernel given by (2.7). If $q = a^2(1+a^2)^{-1}$, then

$$
(2.25) \quad \begin{aligned}
&K_n(2(n-r)+1, x-r+1; 2(n-s)+1, y-s+1) \\
&= \left[\binom{n}{y}\binom{n}{x}^{-1}\right]^{1/2} (-1)^{r-s} L_{n,q}(r, x; s, y),
\end{aligned}
$$

where $K_n$ is the extended Krawtchouk kernel defined by (2.16)–(2.20). The prefactor in (2.25) cancels in all determinants and is not important. Hence, we can just as well think of $L_{n,q}$ as an extended Krawtchouk kernel. We will prove (2.25) in Section 5.

Let us briefly comment on some asymptotic results related to the behavior of the domino tiling of the Aztec diamond close to the point of tangency of the Arctic circle, that is, we consider the positions of the zig-zag particles on the $r$th line where $r$ stays fixed as $n \to \infty$. The corresponding limit for the normalized Krawtchouk polynomials is [19],

$$(2.26) \qquad \lim_{n \to \infty} p_k(qn + \xi\sqrt{2npq}; q, n) = (-1)^k h_k(\xi) \pi^{1/4},$$

where $h_k(x) = 2^{-k/2}(k!)^{-1/2}\pi^{-1/4}H_k(x)$ are the normalized Hermite polynomials. Define an extended Hermite kernel by

$$(2.27) \qquad K_{\mathrm{H,I}}(r,\xi;s,\eta) = \sum_{j=-\infty}^{-1} \sqrt{\frac{(s+j)!}{(r+j)!}} h_{r+j}(\xi) h_{s+j}(\eta) e^{-(\xi^2+\eta^2)/2}$$

if $r \leq s$,

$$(2.28) \qquad K_{\mathrm{H,I}}(r,\xi;s,\eta) = -\sum_{j=0}^{\infty} \sqrt{\frac{(s+j)!}{(r+j)!}} h_{r+j}(\xi) h_{s+j}(\eta) e^{-(\xi^2+\eta^2)/2}$$

if $r > s$. Here $h_k \equiv 0$ if $k < 0$. Note that (2.27) is the ordinary Hermite kernel of $\mathrm{GUE}(r)$ if $r = s$. Taking the limit (2.26) formally in (2.23) and (2.24) gives

$$(2.29) \quad \begin{aligned} &\lim_{n \to \infty} (-n)^{s-r}\sqrt{2npq} L_{n,q}(r, qn + \xi\sqrt{2npq}; s, qn + \eta\sqrt{2npq}) \\ &\qquad = K_{\mathrm{H,I}}(r,\xi;s,\eta). \end{aligned}$$

The kernel $K_{\mathrm{H,I}}$ also occurs in the following problem. Consider a random matrix $A$ from $\mathrm{GUE}(n)$, and let $A_k$ be the $k \times k$ upper left corner of $A$, $1 \leq k \leq n$. The kernel $K_{\mathrm{H,I}}$ then describes the correlations between the eigenvalues $\lambda_j(k,n)$ of $A_k$, $1 \leq j \leq k$, $1 \leq k \leq n$, which form a determinantal process. This follows with some work from [3] and [16]. For the Aztec diamond this leads to the result that the zig-zag particles along the lines $1, \ldots, r$, $r$ fixed, $n \to \infty$, behave like these eigenvalues, see [21]. Note that there is also a different extended Hermite kernel which describes the correlations of the eigenvalues of Dyson's Brownian motion at different times, namely,

$$(2.30) \quad K_{\mathrm{H,II}}(\tau,x;\sigma,y) = \sum_{k=-\infty}^{-1} e^{k(\tau-\sigma)} h_{n+k}(x) h_{n+k}(y) e^{-(x^2+y^2)/2}$$

if $\tau \geq \sigma$, and

$$(2.31) \quad K_{\mathrm{H,II}}(\tau,x;\sigma,y) = -\sum_{k=0}^{\infty} e^{k(\tau-\sigma)} h_{n+k}(x) h_{n+k}(y) e^{-(x^2+y^2)/2}$$

if $\tau < \sigma$. Again, (2.30) reduces to the ordinary Hermite kernel for $\mathrm{GUE}(n)$ if $\tau = \sigma$.

**3. Asymptotics.** The rescaled variables $\xi$, $\xi'$, $\tau$ and $\tau'$ are defined by

$$2r = b_n(\tau) \doteq n(1 + 1/\sqrt{2}) + 2^{-1/6}\tau n^{2/3} \doteq n(\rho + 1),$$
$$2s = b_n(\tau') \doteq n(\rho' + 1),$$
$$x = n(\rho + \sqrt{\rho^2 - 1})/2 + 2^{-5/6}\xi n^{1/3} \doteq n\alpha(\rho) + 2^{-5/6}\xi n^{1/3},$$
$$y = n(\rho' + \sqrt{\rho'^2 - 1})/2 + 2^{-5/6}\xi' n^{1/3}.$$

If we write $a_n(\tau) = n/\sqrt{2} - 2^{-5/6}\tau^2 n^{1/3}$, then $n(\rho + \sqrt{\rho^2 - 1})/2 = a_n(\tau) + \cdots$. Set

$$(3.1) \quad \begin{aligned} K_N^*(2r, x; 2s, y) \\ = (\sqrt{2} - 1)^{x-y+2(s-r)} e^{\xi\tau - \xi'\tau' - (1/3)\tau^3 + (1/3)\tau'^3} K_n(2r, x; 2s, y), \end{aligned}$$

where $K_n$ is given by (2.16). We can use $K_n^*$ as our correlation kernel instead of $K_n$, since all determinants are unchanged. Set $c_n = 2^{-5/6} n^{1/3}$ and define the rescaled kernel

$$(3.2) \quad L_n^*(\tau, \xi; \tau', \xi') = c_n K_n^*(b_n(\tau), a_n(\tau) + c_n \xi; b_n(\tau'), a_n(\tau') + c_n \xi').$$

The basic asymptotic result is contained in the next lemma which will be proved below.

LEMMA 3.1. (a) *Uniformly for $\xi, \xi', \tau, \tau'$ in a compact set,*

$$(3.3) \qquad \lim_{N \to \infty} L_n^*(\tau, \xi; \tau', \xi') = A(\xi, \tau; \xi', \tau'),$$

*where $A$ is the Airy kernel defined by* (1.2).

(b) *Fix $M > 0$. There are positive constants $d_1, d_2, N$, which depend only on $M$ such that if $\tau \geq \tau'$,*

$$(3.4) \qquad |L_n^*(\tau, \xi; \tau', \xi')| \leq d_1 e^{-d_2(\xi + \xi')}$$

*and if $\tau < \tau'$,*

$$(3.5) \qquad |L_n^*(\tau, \xi; \tau', \xi')| \leq \frac{d_1}{\sqrt{\tau' - \tau}} e^{-d_2(\tau' - \tau)(\xi + \xi')}$$

*for all $n \geq N$, $\tau, \tau' \in [-M, M]$ and $\xi, \xi' \in [-M, \infty)$.*

The convergence of Fredholm determinants will be proved using their Fredholm expansions.

LEMMA 3.2. *Let $\tau_1 < \cdots < \tau_m$ and $\gamma_1, \ldots, \gamma_m \in \mathbb{R}$ be given. Set $\Lambda_m = \{\tau_1, \ldots, \tau_m\}$ and let $\nu$ be the counting measure on $\Lambda_m \times \mathbb{Z}$. Also, write*

$g_n(\tau_j, x) = -\mathbb{1}_{(a_n(\tau_j) + \gamma_j c_n, \infty)}(x)$, $x \in \mathbb{Z}$. Then

$$\lim_{n \to \infty} \sum_{k=0}^{n} \frac{1}{k!} \int_{(\Lambda_m \times \mathbb{Z})^k} \det(K_n(z_i, z_j))_{i,j=1}^{k} \prod_{j=1}^{k} g_n(z_j) \, d^k \nu(z) \tag{3.6}$$
$$= \sum_{k=0}^{\infty} \frac{1}{k!} \int_{(\Lambda_m \times \mathbb{Z})^k} \det(A(z_i, z_j))_{i,j=1}^{k} \prod_{j=1}^{k} f(z_j) \, d^k \mu(z),$$

where $\mu = \mu_1 \otimes \mu_2$, $\mu_1$ is a counting measure on $\Lambda_m$, $\mu_2$ is a Lebesgue measure on $\mathbb{R}$ and $f(\tau_j, x) = -\mathbb{1}_{(\gamma_j, \infty)}(x)$, $1 \leq j \leq m$.

PROOF. The proof has many similarities with the proof of Lemma 3.1 in [13], where more details can be found. The sum in the left-hand side of (3.6) can be written

$$\sum_{k=0}^{n} \frac{1}{k!} \int_{(\Lambda_m \times \mathbb{Z})^k} \frac{1}{c_n^k} \det \left( L_n^* \left( \tau(z_i), \frac{x(z_i) - a_n(\tau(z_i))}{c_n}; \right. \right.$$
$$\left. \left. \tau(z_j), \frac{x(z_j) - a_n(\tau(z_j))}{c_n} \right) \right)_{i,j=1}^{k} \tag{3.7}$$
$$\times \prod_{j=1}^{k} f \left( \tau(z_j), \frac{x(z_j) - a_n(\tau(z_j))}{c_n} \right) d^k \nu(z),$$

where $\tau(z_j) = \tau_j$, $x(z_j) = x_j$ if $z_j = (\tau_j, x_j)$. Find $M > 0$ so that $\tau_j \in [-M, M]$, $\gamma_j \geq -M$, $j = 1, \ldots, m$. Then, by (3.4), (3.5) and the Hadamard inequality, the determinant in (3.7) is, $d_3, d_4 > 0$,

$$\leq d_3^k \prod_{i=1}^{k} \left( \sum_{j=1}^{k} \exp \left( -2d_4 \left( \frac{x(z_i) - a_n(\tau(z_i))}{c_n} + \frac{x(z_j) - a_n(\tau(z_j))}{c_n} \right) \right) \right)^{1/2}$$

$$\leq d_3^k e^{kMd_4} k^{k/2} \prod_{i=1}^{k} \exp \left( -d_4 \left( \frac{x(z_i) - a_n(\tau(z_i))}{c_n} \right) \right)$$

for those $z_i$ that contribute to the integral in (3.7). Using this estimate, we see that, given $\varepsilon > 0$, the part of the sum in (3.7) where $k > N$ is

$$\leq \sum_{k=N}^{\infty} \frac{1}{k!} C^k k^{k/2} \left( \frac{1}{c_n} \sum_{j=1}^{m} \sum_{x=a_n(\tau_j) - Mc_n}^{\infty} e^{-d_4(x - a_n(\tau_j))/c_n} \right)^k \leq \varepsilon$$

if $N$ is chosen large enough. Similarly, using the estimates (3.4), (3.5) and the Hadamard inequality again, we can restrict the integration in (3.7) to $\Lambda_m \times ((-\infty, Lc_n] \cap \mathbb{Z})$, with an error $< \varepsilon$ by choosing $L$ large enough. We can then use (3.3) to see that what remains converges to

$$\sum_{k=0}^{N} \frac{1}{k!} \int_{(\Lambda_m \times (-\infty, L])^k} \det(A(z_i, z_j))_{i,j=1}^{k} \prod_{j=1}^{k} f(z_j) \, d^k \mu(z). \tag{3.8}$$

The extended Airy kernel satisfies estimates like (3.4) and (3.5), as can be proved using standard estimates of the Airy function. The same type of argument as above then shows that (3.8) approximates the right-hand side of (3.6), with an error $< \varepsilon$ provided $N$ and $L$ are chosen large enough. This completes the proof. $\square$

We can now give the proof of Theorem 1.1.

PROOF OF THEOREM 1.1. From (2.22) it follows that, with $2r_j = b_n(\tau_j)$,

$$
\begin{aligned}
&\mathbb{P}\left[\max_{1 \leq k \leq N} x_k^{2r_j} \leq a_n(\tau_j) + \gamma_j c_n, 1 \leq j \leq m\right] \\
&\quad = \det(I + gK_n)_{L^2(\Lambda_m \times \mathbb{Z}, \nu)} \\
&\quad = \sum_{k=0}^{n} \frac{1}{k!} \int_{(\Lambda_m \times \mathbb{Z})^k} \det(K_n(z_i, z_j))_{i,j=1}^{k} \prod_{j=1}^{k} g_n(z_j) \, d^k \nu(z),
\end{aligned}
\tag{3.9}
$$

with $g$ as in Lemma 3.2. Now the right-hand side of (3.9) converges to

$$
\sum_{k=0}^{\infty} \frac{1}{k!} \int_{(\Lambda_m \times \mathbb{R})^k} \det(A(z_i, z_j))_{i,j=1}^{k} \prod_{j=1}^{k} f(z_j) \, d^k \mu(z)
= \det(I - f^{1/2} A f^{1/2})_{L^2(\Lambda_m \times \mathbb{R})},
\tag{3.10}
$$

where $f(\tau_j, x) = \mathbb{1}_{(\gamma_j, \infty)}(x)$. Combining (3.9), (3.10) and (1.3), we see that, $N \geq n$,

$$
\begin{aligned}
\lim_{n \to \infty} \mathbb{P}_{N,n}&\left[\max_{1 \leq k \leq N} x_k^{2r_j} \leq n/\sqrt{2} + 2^{-5/6}(\gamma_j - \tau_j^2) n^{1/3}, 1 \leq j \leq m\right] \\
&= \mathbb{P}[\mathcal{A}(\tau_j) \leq \gamma_j, 1 \leq j \leq m].
\end{aligned}
\tag{3.11}
$$

To conclude the proof we use (2.5). Note that in (2.5) the variables are ordered so that $x_1^{2r_j} = \max_{1 \leq k \leq N} x_k^{2r_j}$. Write $x_1^{2r_j} = n/\sqrt{2} + c_n \zeta_j$. Since $2r_j = n(1 + 1/\sqrt{2}) + 2^{-1/6} \tau_j n^{2/3}$, we get

$$
\frac{X_n(2^{-1/6} n^{2/3}(\tau_j - 2^{-2/3} n^{-1/3} \zeta_j)) - n/\sqrt{2} + 1/2}{2^{-5/6} n^{1/3}} = \zeta_j.
\tag{3.12}
$$

Combining (3.11) and (3.12), we obtain (1.4) and the theorem is proved. Note that the variation in the argument of $X_n$ due to $\zeta_j$ is negligible. To be more precise, we could use the uniformity in $\tau_1, \ldots, \tau_m$ in our convegence estimates. $\square$

We still have to prove our basic asymptotic lemma.

PROOF OF LEMMA 3.1. Set

$$
F_{r,x}(z) = z^{x-r}(1-z)^{n-r}(1+z)^r.
$$

Then, by (2.21),

$$(3.13) \quad \widetilde{K}_n(2r,x;2s,y) = \frac{1}{(2\pi i)^2} \int_\Gamma \frac{dz}{z} \int_\gamma \frac{dw}{w} \frac{F_{s,y}(w)}{F_{r,x}(z)} \frac{z}{z-w}$$

and, furthermore,

$$(3.14) \quad \phi_{2r,2s}(x,y) = \frac{1}{2\pi i} \int_\gamma z^{y-x+r-s} \left(\frac{1-z}{1+z}\right)^{r-s} \frac{dz}{z},$$

if $r < s$. Set

$$f(z) = \left(\alpha(\rho) - \frac{1+\rho}{2}\right)\log z + \frac{1-\rho}{2}\log(1-z) + \frac{1+\rho}{2}\log(1+z)$$

so that

$$(3.15) \quad F_{r,x}(z) = e^{2^{-5/6}\xi n^{1/3}\log z + nf(z)}.$$

A straightforward computation shows that $f'(z)$ has a double zero $z_c = z_c(\rho) \doteq \rho(1+\sqrt{1-\rho^2})^{-1}$. It is this double zero condition that specifies $\alpha(\rho)$, and corresponds to the arctic circle. The geometrical considerations and the choice of the place where we want to show convergence to the Airy process leads to $\rho = 1/\sqrt{2} + 2^{-1/6}\tau n^{-1/3}$. (The exact number $2^{-1/6}$ comes out of the computations below.) We get $z_c = \sqrt{2} - 1 + \tau/dn^{1/3} + \cdots$, where $d = 2^{-5/6}(1+\sqrt{2})$. Also, we write $z'_c = z_c(\rho')$.

As our paths of integration we will take

$$z(t) = z_c + \frac{\eta + it}{dn^{1/3}}, \qquad t \in \mathbb{R},$$

$$w(t') = \left(z'_c - \frac{\eta}{dn^{1/3}}\right) e^{it'/d'n^{1/3}}, \qquad |t'| \leq \pi d' n^{1/3},$$

where $d' = d(\sqrt{2}-1) = 2^{-5/6}$, and $\eta > 0$ is such that $\tau - \tau' + 2\eta > 0$.

CLAIM 3.3. *Set*

$$(3.16) \quad \Lambda_n(\eta) = \frac{F_{s,y}(z'_c - \eta/(dn^{1/3}))}{F_{r,x}(z_c + \eta/(dn^{1/3}))}(\sqrt{2}-1)^{x-y+2(s-r)}.$$

*Then with the choices of $r,s,x,y$ as above:*

(i)

$$(3.17) \quad \lim_{n\to\infty} \Lambda_n(\eta) = e^{\xi'\tau' - \xi\tau + \tau^3/3 - \tau'^3/3}$$

*uniformly for $\xi, \xi', \tau, \tau', \eta$ in a compact set, and*

(ii) *for any $\xi, \xi' \geq -M$, $|\tau|, |\tau'| \leq M$, $0 < \eta < 3M$, there is a constant, which only depends on $M$, such that*

$$(3.18) \quad |\Lambda_n(\eta)| \leq Ce^{\xi'\tau' - \xi\tau - \eta(\xi + \xi')}.$$

PROOF. Write
$$w_c = z_c + \eta/dn^{1/3} = \sqrt{2} - 1 + (\tau + \eta)/dn^{1/3} + \cdots,$$
$$w'_c = z'_c - \eta/dn^{1/3} = \sqrt{2} - 1 + (\tau - \eta)/dn^{1/3} + \cdots.$$

The higher-order contributions can be included in $\tau, \tau'$ by slightly changing their values. Then,

$$(\sqrt{2} - 1)^{2r-x} F_{r,x}(w_c)$$
$$= w_c^{x-r}(1 - w_c)^{n-r}(1 + w_c)^r$$
$$= (1 + 2^{5/6}(\tau + \eta)n^{-1/3})^{x-r} 2^{n/2}(1 - 2^{1/3}(\tau + \eta)n^{-1/3})^{n-r}$$
$$\times (1 + (\sqrt{2} - 1)2^{1/3}(\tau + \eta)N^{-1/3})^r.$$

We can now insert $x = n/\sqrt{2} + 2^{-5/6}(\xi - \tau^2)n^{1/3} + \cdots$ and $r = n(1/2 + 1/2\sqrt{2}) + 2^{-7/6}\tau n^{2/3}$ into this expression. Also, we can write an exactly analogous expression for $(\sqrt{2} - 1)^{2s-y} F_{s,y}(w'_c)$. A straightforward, but somewhat lengthy, computation now gives (3.17) and (3.18). □

CLAIM 3.4. *Let $r, s, x, y$ be as above and set*
$$\Omega_n = \{(t, t') \in \mathbb{R} \times [-\pi d' n^{1/3}, \pi d' n^{1/3}]; |t| \geq n^{1/3-\delta}, |t'| \geq n^{1/3-\delta}\}$$

*for some fixed $\delta > 0$. Then, uniformly for $\xi, \xi', \tau, \tau'$ in a compact set,*

$$(3.19) \quad \lim_{n \to \infty} \int_{\Omega_n} \frac{z'(t)w'(t')}{w(t')(z(t) - w(t'))} \frac{F_{s,y}(w(t'))}{F_{r,x}(z(t))})(\sqrt{2} - 1)^{x-y+2(s-r)} \, dt \, dt' = 0.$$

The claim will follow from the estimates used to prove Lemma 3.1(b) below. Let us accept it for the moment. Expanding $f(z)$ around $z_c$, we get

$$(3.20) \quad F_{r,x}(z(t)) = e^{-i\xi(-t+i\eta) - (i/3)(-t+i\eta)^3 + r_n^{(1)}(t)} F_{r,x}(z_c),$$

where $r_n^{(1)}(t) \to 0$ as $n \to \infty$ uniformly for $|t| \leq n^{1/3-\delta}$ and $\xi, \tau$ in a compact set. Similarly,

$$(3.21) \quad F_{s,y}(w(t')) = e^{-i\xi'(t'+i\eta) - (i/3)(t'+i\eta)^3 + r_n^{(2)}(t)} F_{s,y}(z'_c),$$

where $r_n^{(2)}(t) \to 0$ as $n \to \infty$ uniformly for $|t'| \leq n^{1/3-\delta}$ and $\xi', \tau'$ in a compact set. If we insert this into (3.13) and use Claims 3.3 and 3.4, we obtain, after changing $t$ to $-t$,

$$(3.22) \quad \begin{aligned} &\lim_{n \to \infty} (\sqrt{2} - 1)^{x-y+2(s-r)} 2^{-5/6} n^{1/3} \widetilde{K}_n(2r, x; 2s, y) \\ &= -\frac{e^{(1/3)(\tau^3 - \tau'^3) - \xi\tau + \xi'\tau'}}{4\pi^2} \\ &\quad \times \int_{\mathbb{R}} dt \int_{\mathbb{R}} dt' \frac{e^{i\xi(t+i\eta) + i\xi'(t'+i\eta) + i((t+i\eta)^3 + (t'+i\eta)^3)/3}}{\tau' - \tau + i(t + t' + 2\eta)} \\ &= e^{(1/3)(\tau^3 - \tau'^3) - \xi\tau + \xi'\tau'} \widetilde{A}(\tau, \xi; \tau', \xi'), \end{aligned}$$

by Proposition 2.3 in [16]. In order to prove Lemma 3.1(b) and also Claim 3.4, we need some estimates.

Write

(3.23) $$|F_{r,x}(z(t))| = |z(t)|^{2^{-5/6}\xi n^{1/3}}|F_{r,x_0}(z(t))|,$$

where $x_0 = n\alpha(\rho) = n/\sqrt{2} - 2^{-5/6}\tau^2 n^{1/3} + \cdots$. Now,

$$|F_{r,x_0}(\lambda+iu)|^2 = (\lambda^2+u^2)^{x_0-r}((1-\lambda)^2+u^2)^{n-r}((1+\lambda)^2+u^2)^r$$

and, hence,

$$\frac{d}{du}\log|F_{r,x_0}(\lambda+iu)|^2 = 2u\frac{p(u)}{q(u)},$$

where

$$p(u) = (x_0-r)[(1-\lambda)^2+u^2][(1+\lambda)^2+u^2]$$
$$+ (n-r)[\lambda^2+u^2][(1+\lambda)^2+u^2] + r[\lambda^2+u^2][(1-\lambda)^2+u^2]$$

and

$$q(u) = (\lambda^2+u^2)((1+\lambda)^2+u^2)((1-\lambda)^2+u^2).$$

We can now insert $(x_0-r)/n = -1/2 + 1/2\sqrt{2} - 2^{-7/6}\tau n^{-1/3} + \cdots, r/n = 1/2 + 1/2\sqrt{2} + 2^{-7/6}\tau n^{-1/3}$ and $\lambda = \sqrt{2} - 1 + 2^{-5/6}(\sqrt{2}-1)(\tau+\eta)n^{-1/3}$ into $p(u)$. After some computation, we get

$$\frac{1}{n}p(u) = 2^{4/3}(\sqrt{2}-1)^2 \eta n^{-1/3} + (\sqrt{2}-1)u^2 + (1/1/\sqrt{2})u^4,$$

up to negligible contributions. Assume first that $0 \leq u \leq 1$; the case $-1 \leq u \leq 0$ is analogous. Then there is a numerical constant $c_0$ such that

$$\frac{d}{du}\log|F_{r,x_0}(\lambda+iu)| \geq 2c_0\eta n^{2/3}u,$$

and, consequently,

$$|F_{r,x_0}(\lambda+iu)| \geq |F_{r,x_0}(\lambda)|e^{c_0\eta n^{2/3}u^2},$$

for $|u| \leq 1$. (Below $c_0$ denotes a positive numerical constant the value of which may change.) This gives, $|t| \leq dn^{1/3}$,

(3.24) $$|F_{r,x_0}(z(t))| \geq |F_{r,x_0}(z(0))|e^{c_0\eta t^2}.$$

If $u \geq 1$, we use $p(u) \geq nu^4$ and $q(u) \leq c_0 u^6$ and, hence,

$$\frac{d}{du}\log|F_{r,x_0}(\lambda+iu)| \geq c_0 n/u.$$

Combining this with (3.24), we obtain

$$|F_{r,x_0}(z(t))| \geq |F_{r,x_0}(z(0))|e^{c_0\eta n^{2/3}}\left|\frac{t}{dn^{1/3}}\right|^{c_0 n} \quad (3.25)$$

for $|t| \geq dn^{1/3}$. Looking back at (3.23), we see that we must also consider

$$\left|\frac{z(t)}{z(0)}\right|^{2^{-5/6}\xi n^{1/3}} = e^{2^{-1/6}t^2\xi n^{-1/3}+\cdots}. \quad (3.26)$$

If $|\xi| \leq M$, the factor (3.26) can be absorbed into the other factors in the estimates above by slightly changing the constant $c_0$. If $\xi \geq M > 0$, then the expression in (3.26) is $\geq 1$ and we use just this trivial estimate. Combining the estimates, we obtain

$$|F_{r,x}(z(t))| \geq |F_{r,x}(z(0))|g_n(t), \quad (3.27)$$

where

$$g_n(t) = \begin{cases} \exp(c_0\eta t^2), & \text{if } |t| \leq dn^{1/3}, \\ \exp(c_0\eta n^{2/3})|t/dn^{1/3}|^{c_0 n}, & \text{if } |t| \geq dn^{1/3}. \end{cases}$$

We also have to estimate the $w$-integral. Consider

$$|F_{s,y}(w)| = |w|^{y-s}|1-w|^{n-s}|1+w|^s,$$

set $w = \alpha e^{i\theta}$ and define $g(\theta) = (n-s)\log|1-\alpha e^{i\theta}|^2 + s\log|1+\alpha e^{i\theta}|^2$. Then,

$$g'(\theta) = 2\alpha n \sin\theta \frac{(1-2\beta)(1+\alpha^2) + 2\alpha\cos\theta}{(1+\alpha^2)^2 - 4\alpha^2\cos^2\theta}, \quad (3.28)$$

where $2\beta = 2s/n = 1 + 1/\sqrt{2} + 2^{-1/6}\tau' n^{-1/3}$. The numerator in (3.28) is $\leq 0$ for all $\theta$ if $\alpha < z'_c = \sqrt{2} - 1 + 2^{-5/6}(\sqrt{2}-1)\tau' n^{-1/3}$. Hence, this will be satisfied if we choose

$$w(t') = \left(z'_c - \frac{\eta}{dn^{1/3}}\right)e^{it'/d'n^{1/3}},$$

with $\eta > 0$ as above. A computation gives $g(\theta) - g(0) \leq -c_0\theta^2 n^{2/3}\eta$ and, thus,

$$|F_{s,y}(w(t'))| \leq \left|F_{s,y}\left(z'_c - \frac{\eta}{dn^{1/3}}\right)\right|e^{-c_0 t'^2\eta}, \quad (3.29)$$

for $|t'| \leq \pi d' n^{1/3}$. Combining the estimates (3.27) and (3.29), we get

$$|\widetilde{K}_n(2r,x;2s,y)|$$
$$\leq \frac{C}{n^{1/3}}\frac{|F_{s,y}(z'_c - \eta/(dn^{1/3}))|}{|F_{r,x}(z_c + \eta/(dn^{1/3}))|}\left(\int_{-\infty}^{\infty} g_n(t)^{-1}\,dt\right)\left(\int_{-\pi d'n^{1/3}}^{\pi d'n^{1/3}} e^{-c_0 t'^2\eta}\,dt'\right),$$
$$\leq Ce^{\xi'\tau'-\xi\tau-\eta(\xi+\xi')},$$

where we have used (3.18) in the last inequality. Let $\eta = 2\eta_0$. Then, $\xi'(\tau' - \eta) - \xi(\tau + \eta) \leq -\eta_0(\xi + \xi')$, for $\xi, \xi' \geq 0$, if we choose $\eta_0$ so that $\tau' - \eta_0 < 0$, $\tau + \eta_0 > 0$. This proves (3.4) for the $\widetilde{K}$ part, and we also use the same estimates to show Claim 3.4. Notice that $\tau' - \tau - 2\eta < 0$ if $\eta_0$ is chosen large enough.

It remains to consider the $\phi_{r,s}(x,y)$ part of $K_n$ in (2.16). Let $r, x, s, y$ be as above with $\tau' > \tau$. We want to show that, for $\xi, \xi', \tau, \tau'$ in a compact set, we have

$$(3.30) \quad c_n(\sqrt{2}-1)^{x-y+2(s-r)} \phi_{2r,2s}(x,y) \to \frac{1}{\sqrt{4\pi(\tau'-\tau)}} e^{-(\xi'-\xi+\tau^2-\tau'^2)^2/4(\tau'-\tau)}$$

uniformly as $n \to \infty$, and that, given $M > 0$,

$$(3.31) \quad |c_n(\sqrt{2}-1)^{x-y+2(s-r)} e^{\xi\tau-\xi'\tau'} \phi_{2r,2s}(x,y)| \leq \frac{c_0}{\sqrt{\tau'-\tau}} e^{1/4(\tau'-\tau)(\xi+\xi')}$$

for $|\tau|, |\tau'| \leq M$ and $n$ sufficiently large (depending on $\tau' - \tau$).

To prove this we use the integral representation

$$(3.32) \quad \phi_{2r,2s}(x,y) = \frac{1}{2\pi i} \int_{|z|=\alpha} \frac{dz}{z} z^{y-x+r-s} \left(\frac{1+z}{1-z}\right)^{s-r},$$

where $0 < \alpha < 1$. Set $z = \alpha e^{i\theta}$ and $f(\theta) = \log(|1 + \alpha e^{i\theta}|^2 |1 - \alpha e^{i\theta}|^{-2})$. Then

$$f'(\theta) = -4\alpha \sin\theta \frac{1+\alpha^2}{(1+\alpha^2+2\alpha\cos\theta)(1+\alpha^2-2\alpha\cos\theta)},$$

so we have a maximum when $\theta = 0$, and

$$(3.33) \quad f(\theta) - f(0) \leq -c_0 \theta^2$$

for $|\theta| \leq \pi$. We can use this estimate to localize the integral to a small neighborhood of $\theta = 0$ and a standard argument then proves (3.30). From (3.33), we obtain

$$|\phi_{2r,2s}(x,y)| \leq \frac{1}{2\pi} \alpha^{y-x+r-s} \left(\frac{1+\alpha}{1-\alpha}\right)^{s-r} \int_{-\pi}^{\pi} e^{-c_0 n^{2/3}(\tau'-\tau)\theta^2} d\theta$$

$$\leq \alpha^{y-x+r-s} \left(\frac{1+\alpha}{1-\alpha}\right)^{s-r} \frac{1}{n^{1/3}\sqrt{4\pi c_0(\tau'-\tau)}}.$$

Hence,

$$(3.34) \quad \begin{aligned} &c_n(\sqrt{2}-1)^{x-y+2(s-r)} e^{\xi\tau-\xi'\tau'} |\phi_{2r,2s}(x,y)| \\ &\leq \frac{c_0}{\sqrt{\tau'-\tau}} \left(\frac{\alpha}{\sqrt{2}-1}\right)^{y-x+r-s} \left(\frac{1+\alpha}{1-\alpha}(\sqrt{2}-1)\right)^{s-r} e^{\xi\tau-\xi'\tau'}. \end{aligned}$$

Choose $\alpha = (\sqrt{2} - 1)(1 + 2^{-1/6} n^{-1/3}(\tau + \tau'))$. Inserting the expressions for $x, y, r, s$ into (3.34), a computation shows that the right-hand side of (3.34) is

$$\leq \frac{c_0}{\sqrt{\tau' - \tau}} e^{-(\tau' - \tau)(\xi + \xi')/2 + c_1 n^{-1/3}(|\xi| + |\xi'|)}$$

for $|\tau|, |\tau'| \leq M$, where $c_1$ depends on $M$. Hence, if $n$ is large enough, this expression is

$$\leq \frac{c_0}{\sqrt{\tau' - \tau}} e^{-(\tau' - \tau)(\xi + \xi')/4}$$

since $\tau' - \tau > 0$. This completes the proof of Lemma 3.1(b) and, hence, of the whole lemma. $\square$

**4. The center of the Aztec diamond.** It has been conjectured by Propp, see [5, 6, 12], that in the limit $n \to \infty$, a random domino tiling of the Aztec diamond "looks like" a random domino tiling of the whole plane under the Burton–Pemantle measure (the tiling measure with maximal entropy) [4]. In this section we will discuss and outline a proof of a version of this conjecture. The version given below is in terms of particle configurations and not directly in terms of dimers. A domino tiling of the Aztec diamond or the whole plane can equivalently be described as a dimer covering (perfect matching) of a certain graph. Let $V$ be the set of vertices $(1/2 + j, 1/2 + k)$, $j, k \in \mathbb{Z}$, in CS-I and $E$ the set of edges between nearest neighbors. A domino tiling of the whole plane is the same as a dimer covering of $(V, E)$. It will be convenient to write the coordinates as complex numbers. Color the point $1/2 + j + i(k + 1/2)$ black if $j + k$ is even and white otherwise. Let $w$ be a white vertex in $V$ and give an edge between $w$ and $w + z$, where $z = \pm 1, \pm i$, the *weight* $z$. Let $D \subseteq E$ be a finite subset of edges in $E$ and assume that the edges in $D$ cover the black vertices $b_1, \ldots, b_m$ and the white vertices $w_1, \ldots, w_n$. It is proved in [18], using techniques going back to Kasteleyn [17], that under the Burton–Pemantle measure $\mu$, the probability of the event $U_D$ that the edges in $D$ belong to the dimers of a dimer covering of $(V, E)$ equals

$$(4.1) \qquad \mu(U_D) = a_D \det(P(b_j - w_k))_{j,k=1}^m,$$

where $a_d$ is the product of the weights of the edges in $D$ and

$$(4.2) \qquad P(x + iy) = \frac{1}{4\pi^2} \int_{-\pi}^{\pi} \int_{-\pi}^{\pi} \frac{e^{i(x\theta - y\phi)}}{2i \sin \theta + 2 \sin \phi} \, d\theta \, d\phi.$$

Let $A_n$ be the Aztec diamond region as before, set $V_n = A_n \cap V$ and let $E_n$ be the edges between nearest neighbors in $E_n$. Then dimer coverings of the graph $(V_n, E_n)$ are in one-to-one correspondence with domino tilings of $A_n$. Assume that $n$ is odd so that the square with center $1/2 + j + i(k + 1/2)$ is black if and only if $j + k$ is even. The case of even $n$ is analogous. The red and green particles defined in Section 1 lie in $V_n$. It follows from their definition that there is a red particle at $v$ if and only if $v$ is white and there is a dimer from $v - 1$ to $v$ or from $v - i$ to $v$. There is a green particle at $v$ if and only if $v$ is black and there is a dimer from $v$ to $v + 1$ or from $v$ to $v + i$. The discussion in Section 2 also shows that these particles uniquely determine the domino tiling of $A_n$. We can similarly associate a red/green particle process with domino tilings of the whole plane by using the relation to dimers just described. From (4.1) we can then compute the probability of seeing red/green particles at specified points. Note that all particles on the line $y_1 = -x_1 + 2\ell + 1$ have to be green, whereas all the particles on the line $y_1 = -x_1 + 2\ell$ have to be red. Consider, for simplicity, only green particles. We wish to compute the probability of seeing green particles at $v_j = u_j + 1/2 + i(-u_j + 2\ell_j + 1/2)$, $1 \leq j \leq m$. Then all $v_j$ are black and there has to be a dimer between $v_j$ and $v_j + 1$ or between $v_j$ and $v_j + i$. Hence, this probability is given by, compare Theorem 2.10 in [15],

$$
\begin{aligned}
&\sum_{z_j=1 \text{ or } i} z_1 \cdots z_m \det(P(v_j - (v_j + z_k)))_{j,k=1}^m \\
&\sum_{z_j=1 \text{ or } i} z_1 \cdots z_m \sum_{\sigma \in S_m} \text{sgn}(\sigma) \prod_{j=1}^m P(v_j - v_{\sigma(j)} - z_{\sigma(j)}) \\
(4.3) \quad &\sum_{\sigma \in S_m} \text{sgn}(\sigma) \sum_{z_j=1 \text{ or } i} \prod_{j=1}^m z_{\sigma(j)} P(v_j - v_{\sigma(j)} - z_{\sigma(j)}) \\
&\sum_{\sigma \in S_m} \text{sgn}(\sigma) \prod_{j=1}^m \left( \sum_{z_{\sigma(j)}=1 \text{ or } i} P(v_j - v_{\sigma(j)} - z_{\sigma(j)}) \right) \\
&= \det(P(v_j - v_k - 1) + iP(v_j - v_k - i))_{j,k=1}^m.
\end{aligned}
$$

Set, $v \in \mathbb{C}$,

$$(4.4) \qquad R(v) = P(v-1) + iP(v-i).$$

Consider now the probability of seeing green particles at $v_1, \ldots, v_m$ in the Aztec diamond $A_n$, with $n$ odd and sufficiently large to include these points. This probability can be expressed in terms of the extended Krawtchouk kernel. A green point corresponds to a point $Q + (1/2, 0)$ in CS-II where $Q$ is a last particle position. In CS-I this corresponds to a position, compare (2.5), $(2r - n + 1/2 - x^{2r}, x^{2r} - 1/2)$, which gives $2r = 2\ell + n + 1$, $x^{2r} = 2\ell - u + 1$. Hence, in terms of the extended Krawtchouk kernel, (2.16), the probability of finding green particles at $v_1, \ldots, v_m$ is

$$(4.5) \quad \det(K_n(2\ell_j + n + 1, 2\ell - u_j + 1; 2\ell_k + n + 1, 2\ell_k - u_k + 1))_{j,k=1}^m.$$

We can interpret Propp's conjecture as saying that

$$(4.6) \quad \begin{aligned} &\lim_{n \to \infty} \det(K_n(2\ell_j + n + 1, 2\ell - u_j + 1; 2\ell_k + n + 1, 2\ell_k - u_k + 1))_{j,k=1}^m \\ &= \det(R(u_j - u_k + i(u_k - u_j + 2(\ell_j - \ell_k))))_{j,k=1}^m. \end{aligned}$$

The proof of (4.6) consists of two steps. First, we must compute the asymptotics of the kernel $K_n$ for the appropriate values of the arguments corresponding to the center of the Aztec diamond. We will do this using the approach in [22], but we will not give the details, which are of a rather standard saddle-point argument nature. We then have to show that the resulting limiting expression equals the right-hand side of (4.6). These computations will be presented below. The cases of just red particles or a combination of red and green particles can be treated completely analogously.

It follows from (2.16), (2.21), (2.19) and the residue theorem that

$$
(4.7) \quad \begin{aligned} K_n(2r, x; 2s, y) \\ = \frac{1}{(2\pi i)^2} \int_\Gamma \frac{dz}{z} \int_\gamma \frac{dw}{w} \frac{w^y(1-w)^{n-s}(1+1/w)^s}{z^x(1-z)^{n-r}(1+1/z)^r} \frac{z}{z-w}, \end{aligned}
$$

with $\Gamma: t \to \alpha_2 + it$, $t \in \mathbb{R}$, $\gamma: t \to \alpha_1 e^{it}$, $|t| \le \pi$, where $\alpha_2 > \alpha_1 > 0$ if $r \ge s$ and $\alpha_2 < -\alpha_1$, $\alpha_1 > 0$ if $r < s$. Write $n = 2N - 1$, $x_j = 2\ell_j - u_j + 1$, $1 \le j \le m$. Then,

$$
(4.8) \quad \begin{aligned} K_{2N-1}(2(\ell_j + N), x_j; 2(\ell_k + N), x_k) \\ = \frac{1}{(2\pi i)^2} \int_\Gamma \frac{dz}{z} \int_\gamma \frac{dw}{w} \frac{w^{x_k}(1-w)^{N-1-\ell_k}(1+1/w)^{N+\ell_k}}{z^{x_j}(1-z)^{N-1-\ell_j}(1+1/z)^{N+\ell_j}} \frac{z}{z-w}. \end{aligned}
$$

Write $(1-w)^N(1+1/w)^N = w^{-N}(1-w^2)^N = \exp(NF(w))$, where $f(w) = -\log w + \log(1-w^2)$. Then, $f'(w) = 0$ if and only if $w = \pm i$. Assume that $\ell_j > \ell_k$; the other case is similar. We can move $\Gamma$ so that it goes through $\pm i$, $\Gamma_0: t \to it$. We then pick up an extra contribution from the pole at $z = w$ when $w$ is on $\gamma_+: t \to e^{it}$, $-\pi/2 \le t \le \pi/2$. We obtain

$$
(4.9) \quad \begin{aligned} K_{2N-1}(2(\ell_j + N), x_j; 2(\ell_k + N), x_k) \\ = \frac{1}{2\pi i} \int_{\gamma_+} \frac{dw}{w} w^{x_k - x_j}(1-w)^{\ell_j - \ell_k}(1+1/w)^{\ell_k - \ell_j} \\ + \frac{1}{(2\pi i)^2} \int_{\Gamma_0} \frac{dz}{z} \int_\gamma \frac{dw}{w} \frac{w^{x_k}(1-w)^{N-1-\ell_k}(1+1/w)^{N+\ell_k}}{z^{x_j}(1-z)^{N-1-\ell_j}(1+1/z)^{N+\ell_j}} \frac{z}{z-w}. \end{aligned}
$$

The second integral goes to zero as $N \to \infty$, compare [22], and we obtain

$$
(4.10) \quad \begin{aligned} \lim_{N \to \infty} K_{2N-1}(2(\ell_j + N), 2\ell_j - u_j + 1; 2(\ell_k + N), 2\ell_k - u_k + 1) \\ = \frac{1}{2\pi i} \int_{\gamma_+} \frac{dw}{w} w^{u_j - u_k} \left( \frac{1-w}{w(1+w)} \right)^{\ell_j - \ell_k}. \end{aligned}
$$

It follows from (4.4) and (4.10) that in order to prove (4.6) it suffices to show that

$$
(4.11) \quad \begin{aligned} \frac{1}{2\pi i} \int_{\gamma_+} \frac{dw}{w} w^{u-y} \left( \frac{1-w}{w(1+w)} \right)^{\ell-m} \\ = i^{u-v+2(\ell-m)}[P(u - v - 1 + i(v - u + 2(\ell - m))) \\ + iP(u - v + i(v - u - 1 + 2(\ell - m)))]. \end{aligned}
$$

We will use the following formula. If $x + y \geq 1$ and $x + y$ is odd, then

$$(4.12) \qquad P(x+iy) = \frac{i^{x-1}}{2\pi i} \int_{\gamma_+} w^{(x-y+1)/2} \frac{(w-1)^{(x+y-1)/2}}{(w+1)^{(x+y+1)/2}} \frac{dw}{w}.$$

Assume that we have proved (4.12). If we set $x = u - v - 1$, $y = v - u + 2(\ell - m)$, then the right-hand side of (4.11) is

$$-(-1)^{(x+y+1)/2} i^{x+1} \frac{1}{2\pi i} \int_{\gamma_+} (i^{x-1} + i^{x+1}w) w^{(x-y+1)/2} \frac{(w-1)^{(x+y-1)/2}}{(w+1)^{(x+y+1)/2}} \frac{dw}{w}$$

$$= \frac{1}{2\pi i} \int_{\gamma_+} w^{(x-y+1)/2} \left(\frac{w-1}{w+1}\right)^{(x+y+1)/2} \frac{dw}{w},$$

which is exactly the left-hand side of (4.11). It remains to prove (4.12).

Make the shift $\theta \to \theta - \phi$ in (4.2) to get

$$(4.13) \quad P(x+iy) = \frac{1}{4\pi^2} \int_{-\pi}^{\pi} \left[ \int_{-\pi}^{\pi} \frac{e^{-i(x+y)\phi}}{2i\sin(\theta-\phi) + 2\sin\phi} d\phi \right] e^{ix\theta} d\theta.$$

Since $2i\sin(\theta-\phi) + 2\sin\phi = [e^{i\theta} + i - (e^{-i\theta} + i)e^{2i\phi}]e^{-i\phi}$, we obtain

$$P(x+iy) = \frac{1}{4\pi^2} \int_{-\pi}^{\pi} \left[ \int_{-\pi}^{\pi} \frac{e^{i((x+y-1)/2)2\phi}}{e^{i\theta} + i - (e^{-i\theta} + i)e^{2i\phi}} d\phi \right] e^{ix\theta} d\theta$$

$$= \frac{1}{2\pi} \int_{-\pi}^{\pi} \frac{e^{ix\theta}}{e^{i\theta} + i} \left[ \frac{1}{2\pi i} \int_{|z|=1} \frac{z^{-(x+y-1)/2-1}}{1 - az} dz \right] d\theta,$$

where $a = (e^{-i\theta} + i)(e^{i\theta} + i)^{-1}$. If we make use of the integral, compare [22],

$$\frac{1}{2\pi i} \int_{|z|=1} \frac{z^{-k-1}}{1 - az} dz = \begin{cases} a^k, & \text{if } k \geq 0, |a| < 1, \\ -a^k, & \text{if } k < 0, |a| > 1, \\ 0, & \text{otherwise,} \end{cases}$$

we find

$$P(x+iy) = \frac{1}{2\pi} \int_0^{\pi} \frac{e^{ix\theta}}{e^{i\theta} + i} \left(\frac{e^{-i\theta} + i}{e^{i\theta} + i}\right)^{(x+y-1)/2} d\theta,$$

which is easily seen to equal the right-hand side of (4.12). A similar computation should relate the Kasteleyn kernel for the finite Aztec diamond, computed in [11], to the extended Krawtchouk kernel $K_n$.

**5. Extended Krawtchouk kernel and Krawtchouk polynomials.** In this section we will prove (2.25), which expresses the extended Krawtchouk kernel in terms of Krawtchouk polynomials. We have the formula

$$(5.1) \quad \begin{aligned} & K_n(2r - \varepsilon_1, x; 2s - \varepsilon_2, y) \\ & = \frac{1}{(2\pi i)^2} \int_{\gamma_{r_1}} \frac{dz}{z} \int_{\gamma_{r_2}} \frac{dw}{w} \frac{w^y}{z^x} \frac{(1-aw)^{n-s+\varepsilon_2}(1+a/w)^s}{(1-az)^{n-r+\varepsilon_1}(1+a/z)^r} \frac{z}{z-w}, \end{aligned}$$

where $\varepsilon_1, \varepsilon_2 \in \{0, 1\}$, $a < r_1 < 1/a$ and $0 < r_2, r_1$ if $2r - \varepsilon_1 \geq 2s - \varepsilon_2$, $0 < r_1 < r_2$ if $2r - \varepsilon_1 < 2s - \varepsilon_2$. Here $\gamma_r$ denotes a circle of radius $r$ around the origin. This follows from (2.16)–(2.20) and the residue theorem.

Consider first the case $r \leq s$. Then

$$
(5.2) \quad \begin{aligned} &K_n(2(n-r)+1, x-r+1; 2(n-s)+1, y-s+1) \\ &= \frac{(-1)^{s-r}}{(2\pi i)^2} \int_{\gamma_{r_1}} \frac{dz}{z} \int_{\gamma_{r_2}} \frac{dw}{w} \frac{z^{n-x}}{w^{n-y}} \frac{(aw-1)^s (w+a)^{n-s+1}}{(az-1)^r (z+a)^{n-r+1}} \frac{z}{z-w}. \end{aligned}
$$

Since $|w| < |z|$, we have

$$
\sum_{k=0}^{\infty} \frac{((w+a)/(aw-1))^k}{((z+a)/(az-1))^k} = \frac{(aw-1)(z+a)}{(a^2+1)(w-z)}
$$

and, hence, the integral in the right-hand side of (5.2) can be written

$$
(5.3) \quad \begin{aligned} &-\frac{(-1)^{s-r}(a^2+1)}{(2\pi i)^2} \sum_{k=0}^{\infty} \left( \int_{\gamma_{r_1}} \frac{z^{n-x}}{(az-1)^{r-k}(z+a)^{n-r+k+2}} \frac{dz}{z} \right) \\ &\qquad\qquad\times \left( \int_{\gamma_{r_2}} \frac{(aw-1)^{s-k-1}(w+a)^{n+k-s+1}}{w^{n-y}} \frac{dw}{w} \right). \end{aligned}
$$

In the $z$-integral we make the change of variables $z = (1 - a\zeta)(\zeta + a)^{-1}$, $\zeta = (1 - az)(z + a)^{-1}$. This maps $\gamma_{r_1}$ to a circle surrounding the origin but with the opposite orientation. We can deform this circle to $\gamma_{r_1}$ using Cauchy's theorem. Hence, (5.3) equals

$$
(5.4) \quad \begin{aligned} &\frac{1}{(2\pi i)^2} \sum_{k=0}^{r-1} \frac{1}{(a^2+1)^n} \left( \int_{\gamma_{r_1}} \frac{(1-a\zeta)^{n-x}(\zeta+a)^x}{\zeta^{r-k-1}} \frac{d\zeta}{\zeta} \right) \\ &\qquad\qquad\times \left( \int_{\gamma_{r_2}} \frac{(1-a\omega)^{s-k-1}(\omega+a)^{n+k-s+1}}{\omega^{n-y}} \frac{d\omega}{\omega} \right). \end{aligned}
$$

Here we have also used the fact that, by Cauchy's theorem, the $\zeta$-integral vanishes if $k \geq r$. We now set $\omega = w/a$, $\zeta = a(1+a^2)^{-1}z$ and use $p = (1+a^2)^{-1}$, $q = a^2(1+a^2)^{-1}$. A computation, where we also replace $k$ by $r-1-k$, shows that (5.4) equals

$$
(5.5) \quad \begin{aligned} &\frac{a^{s-r}(1+a^2)^{r-s}a^{x-y}}{(2\pi i)^2} \sum_{k=0}^{r-1} \left( \int_{\gamma_{r_1}} \frac{(1-qz)^{n-x}(1+pz)^x}{z^k} \frac{dz}{z} \right) \\ &\qquad\qquad\times \left( \int_{\gamma_{r_2}} \frac{(1-w)^{s-r+k}(pw+q)^{n+r-s-k}}{w^{n-y}} \frac{dw}{w} \right). \end{aligned}
$$

By Cauchy's integral formula, this is the coefficient of $w^{n-y}$ in

$$
\left( \frac{a}{1+a^2} \right)^{s-r} a^{x-y} \sum_{k=0}^{r-1} (1-w)^{s-r+k}(pw+q)^{n+r-s-k} \frac{1}{2\pi i} \\ \times \int_{\gamma_{r_1}} \frac{(1-qz)^{n-x}(1+pz)^x}{z^k} \frac{dz}{z}.
$$

This can be written

$$
\begin{aligned}
&\left(\frac{a}{1+a^2}\right)^{s-r} a^{x-y} \sum_{k=s-r}^{s-1} (1-w)^k (pw+q)^{n-k} \frac{1}{2\pi i} \\
&\qquad \times \int_{\gamma_{r_1}} \frac{(1-qz)^{n-x}(1+pz)^x}{z^{k+r-s}} \frac{dz}{z} \\
&= \left(\frac{a}{1+a^2}\right)^{s-r} a^{x-y} \sum_{k=s-r}^{s-1} \sum_{y=0}^{n} \delta_{yk} (1-w)^k (pw+q)^{n-k} \frac{1}{2\pi i} \\
&\qquad \times \int_{\gamma_{r_1}} \frac{(1-qz)^{n-x}(1+pz)^x}{z^{k+r-s}} \frac{dz}{z} \\
&= \left(\frac{a}{1+a^2}\right)^{s-r} a^{x-y} \sum_{k=s-r}^{s-1} \binom{n}{k}^{-1} (pq)^{-k} \frac{1}{2\pi i} \int_{\gamma_{r_2}} \frac{d\omega}{\omega^{k+1}} \\
&\qquad \times \sum_{y=0}^{n} \binom{n}{y} (1-w)^y (pq)^y \omega^y (pw+q)^{n-y} \frac{1}{2\pi i} \\
&\qquad \times \int_{\gamma_{r_1}} \frac{(1-qz)^{n-x}(1+pz)^x}{z^{k+r-s}} \frac{dz}{z}.
\end{aligned}
\tag{5.6}
$$

By the binomial theorem, the $y$-sum equals

$$(pq\omega - pqw\omega + pw + q)^n = (pq\omega + q + w(p - pq\omega))^n$$
$$= \sum_{y=0}^{n} \binom{n}{y} (pq\omega + q)^y (p - pq\omega)^{n-y} w^{n-y}.$$

Hence, the coefficient of $w^{n-y}$ in (5.6) equals, by (2.8),

$$
\begin{aligned}
&\sum_{k=s-r}^{s-1} \binom{n}{k}^{-1} (pq)^{-k} w_q(y) \left(\frac{a}{1+a^2}\right)^{s-r} a^{x-y} \\
&\qquad \times (-1)^k (pq)^{k/2} \binom{n}{k}^{1/2} p_k(y;q,n) (-1)^{k+r-s} (pq)^{(k+r-s)/2} \\
&\qquad \times \binom{n}{k+r-s}^{1/2} p_{k+r-s}(x;q,n) \\
&= \left[\binom{n}{y}\binom{n}{x}^{-1}\right]^{1/2} \\
&\qquad \times \sum_{k=s-r}^{s-1} \binom{n}{k}^{-1/2} \binom{n}{k+r-s}^{1/2} p_{k+r-s}(x;q,n) p_k(y;q,n),
\end{aligned}
$$

which proves (2.25) when $r \le s$.

The case when $r > s$ will be handled by reducing to the previous case. In the integral in

$$K_n(2(n-r)+1, x-r+1; 2(n-s)+1, y-s+1)$$
$$= \frac{(-1)^{s-r}}{(2\pi i)^2} \int_{\gamma_{r_1}} \frac{dz}{z} \int_{\gamma_{r_2}} \frac{dw}{w} \frac{z^{n-x}}{w^{n-y}} \frac{(aw-1)^s(w+a)^{n-s+1}}{(az-1)^r(z+a)^{n-r+1}} \frac{z}{z-w},$$

with $a < r_1 < 1/a$, $r_2 > r_1$, we make the change of variables $w \to -1/w$, $z \to -1/z$. If $r_3 = 1/r_1$, $r_4 = 1/r_2$, we see that

$$K_n(2(n-r)+1, x-r+1; 2(n-s)+1, y-s+1)$$
$$= -\frac{(-1)^{y-x+r-s}}{(2\pi i)^2} \int_{\gamma_{r_3}} \frac{dz}{z} \int_{\gamma_{r_4}} \frac{dw}{w} \frac{z^{n-(n-x)}}{w^{n-(n-y)}}$$
$$\times \frac{(1-aw)^{n-s+1}(w+a)^{n-(n-s+1)+1}}{(1-az)^{n-r+1}(z+a)^{n-(n-r+1)+1}} \frac{z}{z-w}$$
$$= -(-1)^{y-x+r-s} K(2(n-(n-r+1))+1, n-x-(n-r+1)+1;$$
$$2(n-(n-s+1))+1, n-y-(n-s+1)+1).$$

By our previous computation this equals

$$-(-1)^{y-x+r-s}\left[\binom{n}{y}\binom{n}{x}^{-1}\right]^{1/2}(-1)^{s-r}\sum_{k=r-s}^{n-s}\binom{n}{k}^{-1/2}\binom{n}{k+s-r}^{1/2}$$
(5.7)
$$\times p_{k+s-r}(n-x; q, n)$$
$$\times p_k(n-y; q, n)$$
$$\times [w_q(n-x)w_q(n-y)]^{1/2}.$$

Using the integral formula, (2.8), for the normalized Krawtchouk polynomials, it is not difficult to show that

(5.8) $\quad p_k(n-x; q, n) = (-1)^{k-x} p^{n/2-x} q^{-n/2+x} p_{n-k}(x; q, n).$

If we use this formula in (5.7) we obtain, after some simplification,

$$-\left[\binom{n}{y}\binom{n}{x}^{-1}\right]^{1/2}(-1)^{r-s}\sum_{k=r-s}^{n-s}\binom{n}{k}^{-1/2}\binom{n}{k+s-r}^{1/2}$$
$$\times p_{n-k+r-s}(x; q, n)$$
$$\times p_{n-k}(y; q, n)[w_q(x)w_q(y)]^{1/2}$$
$$-\left[\binom{n}{y}\binom{n}{x}^{-1}\right]^{1/2}(-1)^{r-s}\sum_{j=0}^{\infty}\binom{n}{j+s}^{-1/2}\binom{n}{j+r}^{1/2}$$
$$\times p_{j+r}(x; q, n)$$
$$\times p_{j+s}(y; q, n)[w_q(x)w_q(y)]^{1/2},$$

where we have put $j = n - k - s$ and extended the summation using $p_k \equiv 0$ if $k < 0$ or $k > n$. This is what we wanted to prove.

**Acknowledgments.** I thank A. Borodin for a helpful discussion concerning extended kernels and E. Nordenstam for help with the pictures.

DEPARTMENT OF MATHEMATICS
ROYAL INSTITUTE OF TECHNOLOGY
S-100 44 STOCKHOLM
SWEDEN
E-MAIL: kurtj@math.kth.se